\documentclass[pdflatex,sn-mathphys-num]{sn-jnl}% Math and Physical Sciences Numbered Reference Style
\usepackage{graphicx}%
\usepackage{multirow}%
\usepackage{amsmath,amssymb,amsfonts}%
\usepackage{amsthm}%
\usepackage{mathrsfs}%
\usepackage[title]{appendix}%
\usepackage{xcolor}%
\usepackage{textcomp}%
\usepackage{manyfoot}%
\usepackage{booktabs}%
\usepackage{algorithm}%
\usepackage{algorithmicx}%
\usepackage{algpseudocode}%
\usepackage{listings}%
\newcommand{\N}{{\mathbb{N}}}

\newcommand{\E}{{\mathbb{E}}}

\newcommand{\R}{{\mathbb{R}}}

\newcommand{\Ln}{{\mathcal{L}_{n}}}

\newcommand{\LnN}{{\mathcal{L}_{n}^{\mathbb{N}}}}

\newcommand{\BV}{{BV_{0}(0,1)}}
\newcommand{\BVM}{{BV_{0,M}(0,1)}}

\newcommand{\cL}{{\mathcal{L}}}

\newcommand{\cU}{{\mathcal{U}}}

\newcommand{\cD}{{\mathcal{D}}}
\newcommand{\cS}{{\mathcal{S}}}
\newcommand{\cW}{{\mathcal{W}}}
\newcommand{\cN}{{\mathcal{N}}}
\newcommand{\cI}{{\mathcal{I}}}
\newcommand{\cQ}{{\mathcal{Q}}}
\newcommand{\cH}{{\mathcal{H}}}
\newcommand{\cF}{{\mathcal{F}}}
\newcommand{\cE}{{\mathcal{E}}}

\newcommand{\cK}{{\mathcal{K}}}
\newcommand{\cB}{{\mathcal{B}}}
\newcommand{\cT}{{\mathcal{T}}}

\newcommand{\cA}{{\mathcal{A}}}
\newcommand{\cC}{{\mathcal{C}}}

\newcommand{\rI}{{\mathrm{I}}}

\newcommand{\id}{{\mathrm{i_{d}}}}

\newcommand{\rc}{{\mathrm{c}}}

\newcommand{\lra}{{\longrightarrow}}
\newcommand{\ra}{{\rightarrow}}

\theoremstyle{thmstyleone}

\newtheorem{rem}{Remark}
\newtheorem{prop}{Proposition}
\newtheorem{corr}{Corollary}
\newtheorem{lem}{Lemma}

\newtheorem{theo}{Theorem}
\newtheorem{defi}{Definition}

\begin{document}

\title[Article Title]{Functional Erd\H{o}s-R\'enyi laws for L\'evy processes}


%\author{\fnm{Dimbihery} \sur{Rabenoro}}
\author*[]{\fnm{Dimbihery} \sur{Rabenoro}}\email{dimbihery.rabenoro@inria.fr}
%\email{dimbihery.rabenoro@inria.fr}



\affil{\orgdiv{LPSM}, \orgname{Sorbonne University}, \orgaddress{\city{Paris}, \country{France}}}

%\affil[2]{\orgdiv{Department}, \orgname{Organization}, \orgaddress{\street{Street}, \city{City}, \postcode{10587}, \state{State}, \country{Country}}}

%\affil[3]{\orgdiv{Department}, \orgname{Organization}, \orgaddress{\street{Street}, \city{City}, \postcode{610101}, \state{State}, \country{Country}}}


\abstract{In this paper we establish functional Erd\H{o}s-Renyi laws for L\'evy processes, i.e. limit theorems for sets of functions on [0,1] associated to their increments. First, we determine precise conditions under which, in a general framework, such a convergence is derived from a large deviations principle for probability measures induced by the sample paths of such a process. Then, by checking that these conditions are fulfilled, we obtain, under two usual assumptions on exponential moments, such limit theorems from well-known large deviations principles.}





\maketitle

\section{Introduction}\label{sect1}

\noindent
For a sequence $(X_{i})_{i \geq 1}$ of i.i.d. $\R$-valued random variables, consider the assumptions:
\begin{center}
$(\cC_{X})$ : $\phi_{X}(\theta) < \infty$ for all $\theta \in \R$ ~;~ $(\cA_{X})$ : $\phi_{X}(\theta) < \infty$ for $\theta$ in a neighborhood of $0$. 
\end{center}

\noindent
where $\phi_{X}(\theta):=\E[\exp(\theta X_{1})]$. The celebrated Erd\H{o}s-Renyi Theorem provides, in \cite{erdos_new_1970}, a law of large numbers for the supremum of moving averages of the $(X_{i})$, under $(A_{X})$. This fundamental result has been extended to so-called functional Erd\H{o}s-Renyi laws [FERL]'s, which are functional limit theorems for the increments of the partial sum process, notably in \cite{borovkov_1990} and \cite{sanchis_functional_1994} under $(C_{X})$, and in \cite{deheuvels1991functional} under $(A_{X})$. When $X_{1}$ has only a semiexponential distribution, a FERL is derived in \cite{gantert_functional_1998}. More recently, in \cite{houdre_erdos-renyi-type_2016} and \cite{kifer_functional_2017}, FERL's are obtained for processes which are variants of the partial sum one.      

\noindent\\
In this paper, we establish a functional limit of this type for L\'evy processes. Given a L\'evy process $Z:=\left\{ Z(t) : t \geq 0 \right\}$, 
for $x \geq 0$ and $\lambda>0$, we consider the standardized increment function $\eta_{x,\lambda}$ of $Z$ defined by  
\begin{equation*}
\eta_{x,\lambda} : [0,1] \lra \R \quad , \quad \eta_{x,\lambda}(s) = \lambda^{-1}(Z(x+s\lambda)-Z(x)). 
\end{equation*}

\noindent
For $c>0$ and $n>1$, let $A_{n}$ be the integer part of $c\log n$. Then, consider the following sets of increment functions:
\begin{equation*}
\Ln:=\left\{\eta_{x,A_{n}} : x \in [0; n-A_{n}] \right\} \quad \textrm{and} \quad \LnN:=\left\{\eta_{x,A_{n}} : x \in [0; n-A_{n}] \cap \N \right\}. 
\end{equation*}

\noindent 
By definition, all increment functions of $Z$ lie a.s. in the Skorokhod space $D(0,1)$ of c\`adl\`ag functions on $[0,1]$. However, we consider the general case of a vector space $\cE$ of functions on $[0,1]$ such that a.s., all increment functions of $Z$ lie in $\cE$. Given a distance $d$ on $\cE$, let $\cD$ be the Hausdorff distance between subsets of $\cE$ induced by $d$: See subsection \ref{PreliminariesMain}. We say that a functional Erd\H{o}s-Renyi law [FERL] for a subsequence of $(\Ln)$ holds in $(\cE, d)$ with limit set $\cK \subset \cE$ when, 
outside a negligible set, the sequence $(\Ln)$ converges, with respect to $\cD$, to $\cK$. 

As expected, such a FERL should derive from a large deviation principle [LDP]. Thus, we consider, for $\lambda>0$, the distribution $P_{\lambda}$ of $\eta_{0, \lambda}=\lambda^{-1}Z(\lambda \cdot)$. Then, well-known LDP's for the family $(P_{\lambda})$ are available, under the following usual assumptions on the moment generating function $\Phi$ of $Z(1)$: 
\begin{center}
$(\cC)$ : $\Phi(\theta) < \infty$ for all $\theta \in \R$ \quad ; \quad $(\cA)$ : $\Phi(\theta) < \infty$ for $\theta$ in a neighborhood of $0$.
\end{center} 

\noindent
Under $(\cC)$, a LDP for $(P_{\lambda})$ in $(D(0,1), \cS)$ is obtained in \cite{varadhan_asymptotic_1966}, where $\cS$ is the Skorokhod topology. When only $(\cA)$ holds, the rate function of this LDP loses a key property: Its sublevel sets are not compact anymore, which is discussed in Remark 2 of \cite{borovkov_1990}. In fact, under additional conditions on $Z$, its increment functions lie a.s. in the space $\BV$ of functions with bounded variations on $[0,1]$ that vanish at $0$. Then, under $(\cA)$, a LDP in $(\BV, \cW)$ is established in \cite{lynch_large_1987}, where $\cW$ is the weak topology. 

\noindent\\
The novelty of this paper is to consider a unified framework, instead of deriving FERL's separately from each of these LDP's. Thus, given any LDP for $(P_{\lambda})$ in $(\cE, \cT)$, where $\cT$ is a topology on $\cE$, we determine a synthetic set of conditions under which a FERL for $(\Ln)$ in $(\cE, d)$ is derived from such a LDP: See subsection  \ref{subsecGenFrame}. This approach is similar to that of \cite{deheuvels_strassen-type_1993}, where the authors characterize the topologies under which the functional law of the iterated logarithm is valid. We now describe these conditions. We first prove that a FERL for $(\Ln)$ holds in $(\cE, d)$ when a FERL for $\left(\LnN \right)$ holds in $(\cE, d)$ and the following condition $(\Delta^{d})$ is fulfilled. 
\begin{center}
$(\Delta^{d})$ : For all $\epsilon>0$, $\left\{ \Delta_{n}^{d}(\epsilon) \textrm{ infinitely often in n } \right\}$ is negligible, 
\end{center}

\noindent
where $\Delta_{n}^{d}(\epsilon) := \left\{ \exists x \in [0, n-A_{n}] \textrm{ ~such that~ } d\left( \eta_{x,A_{n}}, \eta_{\lfloor x \rfloor,A_{n}} \right) \geq \epsilon \right\}$. Then, we identify conditions under which a FERL for $\left(\LnN \right)$ holds. First, it is required that the LDP has right properties, notably that its rate function is good (i.e. its sublevel sets are compact) and convex. Then, other conditions relate to relationships between $\cT$ and the topology induced by $d$. In fact, when $d$ metricizes $\cT$, these conditions for $\left(\LnN \right)$ are reduced to the LDP properties. 

\noindent\\
Then, by checking that these conditions hold, we apply the general result of Theorem \ref{generalTheoN} to derive FERL's from the LDP's of \cite{varadhan_asymptotic_1966} and \cite{lynch_large_1987}. Thus, we obtain in Theorem \ref{ferlPartialCA} FERL's for $\left(\LnN \right)$ under $(\cC)$ and $(\cA)$, respectively in $(D(0,1), d_{\cU})$ and $(\BV, d_{\cH})$, where $d_{\cU}$ is the uniform distance on $D(0,1)$ and $d_{\cH}$ a distance on $\BV$ introduced in \cite{hognas1977characterization}: See $(\ref{aa1})$ in subsection \ref{BVproperties}.

Under $(\cC)$, we also derive in Theorem \ref{ferlFullC} a FERL for the whole $(\Ln)$ in $(D(0,1), d_{\cU})$, by proving additionally that condition $\left( \Delta^{d_{\cU}} \right)$ holds. In fact, a FERL for $\left(\LnN \right)$ holds in $(D(0,1), d_{\cS})$ but condition $\left( \Delta^{d_{\cS}} \right)$ is difficult to check, since $d_{\cS}$ is not easily handleable. On the other hand, the FERL for $\left(\LnN \right)$ in $(\BV, d_{\cH})$ would extend to a FERL for the whole $(\Ln)$ if condition $\left( \Delta^{d_{\cH}} \right)$ holds under $(\cA)$, but this is still an open question. However, it is established in \cite{frolov_unified_2005} that, under $(\cA)$, there exists a constant $\gamma_{c}$ such that  
\begin{equation}\label{ERofFrolov}
\sup\limits_{0 \leq x \leq n - A_{n}} A_{n}^{-1} \left(Z(x + A_{n}) - Z(x) \right) ~{\underset{n \ra \infty}\lra} ~\gamma_{c} \quad \textrm{a.s.}
\end{equation}

\noindent
If a FERL for the whole $(\Ln)$ holds under $(\cA)$, then it would imply the convergence of $(\ref{ERofFrolov})$. Thus, it is not unreasonable to conjecture that such a FERL holds under, at least, a condition weaker than $(\cC)$. 

\noindent\\
Finally, we mention that for other types of increment functions of $Z$, a functional limit theorem [FLT] is derived directly from an invariance principle stated in \cite{hoffmann-jorgensen_random_1994} which provides, under $(\cA)$, an approximation of $Z$ by a Wiener process $W$. Indeed, combining it with a FLT for some increment functions of $W$ obtained in \cite{rvsz_generalization_1979}, we deduce a FLT for increment functions of $Z$ of the form $\nu_{x,\alpha_{n}} : [0,1] \lra \R$ and $\nu_{x,\alpha_{n}}(s) = \beta_{n}^{-1} \left(Z(x+s\alpha_{n})-Z(x) \right)$, where $x \geq 0$, $(\alpha_{n})$ is a sequence of increments larger than $\log n$, i.e. $\alpha_{n} / \log n \ra \infty$ as $n \ra \infty$ and $(\beta_{n})$ a suitable normalizing sequence. Since the increments $A_{n}$ of $\eta_{x,A_{n}}$ are of order $\log n$, one checks readily that this method by invariance principle is not applicable to derive a FLT for the $(\eta_{x,A_{n}})$, which motivates our approach by LDP's.

\noindent\\
This paper is organized as follows. In Sections \ref{sec:2} and \ref{sec:3}, we present preliminaries on L\'evy processes and functional spaces. Our main results are stated in Section \ref{sec:4}, with proofs and technical results detailed in Section \ref{sec:5}.

\section{L\'evy Processes}\label{sec:2}

\noindent
Here, stochastic processes are families of $\R$-valued random variables, indexed by $[0, \infty)$.

\subsection{Basic properties of L\'evy processes}

\begin{defi}
For any interval $\rI$, $D(\rI)$ is the set of c\`adl\`ag functions on $\rI$. On $D([0,1])$, denoted by $D(0,1)$, the Skohorod topology $\mathcal{S}$ is induced by the distance 
\begin{equation}\label{defsko}
d_{\mathcal{S}}(f,g) = \inf\limits_{\nu \in \Lambda} \Big\{ \max \left( \left\| \nu - I \right\|  ;
\left\| f - g\circ \nu \right\| \right) \Big\},
\end{equation}

\noindent
where $\Lambda$ is the class of strictly increasing, continuous mappings of $[0,1]$ onto itself. 

\noindent
A L\'evy process is a stochastic process $Z=\left\{ Z(t) : t \geq 0 \right\}$ such that: $(i)$ $Z(0)=0$ a.s. $(ii)$ $Z$ has stationary and independent increments and is stochastically continuous. $(iii)$ Its sample paths lie in $D([0,\infty))$ a.s.
\end{defi}

\noindent
For any L\'evy process $Z$, the distribution of $Z(1)$ is an infinitely divisible distribution [IDD]. We define hereafter a probability space on which L\'evy processes are defined. We need the following preliminary definition. 

\begin{defi}
For all interval $\rI$, the class $\mathfrak{C}_{\rI}$ of cylinder sets is formed by sets of $\R^{\rI}$ of the form $\bigcap\limits_{i=1}^{n} \left\{ f \in \R^{\rI} : f(t_{i}) \in B_{i} \right\}$, for $n \geq 1$, $(t_{i})_{1 \leq i \leq n} \in \rI^{n}$ and $(B_{i})_{1 \leq i \leq n} \in (\cB_{\R})^{n}$, where $\cB_{\R}$ is the Borel $\sigma$-algebra of $\R$. 
\end{defi}

\noindent
Let $\Omega = \R^{[0,\infty)}$ and $\cF$ the $\sigma$-algebra on $\Omega$ generated by $\mathfrak{C}_{[0,\infty)}$. Given an IDD $\mu$ on $\R$, there exists a unique probability measure $P$ on $(\Omega, \cF)$ and a L\'evy process $Z$ on $(\Omega, \cF, P)$ such that $Z(1)$ is distributed as $\mu$: See \cite{sato_levy_2013}. Finally, we recall that a L\'evy process is fully detemined by its \textit{generating triplet} $(A, \nu, \gamma)$, where $A \geq 0$ is the variance of its Gaussian component, $\nu$ is a measure on $\R$ such that $\nu(\left\{ 0 \right\})=0$ and $\int_{\R} \left(x^{2} \wedge 1 \right) \nu(dx) < \infty$, and $\gamma \in \R$.

\begin{lem}\label{pathsBV}
By Theorem 21.9 in \cite{sato_levy_2013}, if $Z$ is a L\'evy process whose generating triplet $\left(A, \nu, \gamma \right)$ satisfies that $A=0$ and $\int_{|x| \leq 1} |x|\nu(dx) < \infty$, then, a.s., for any $t \in (0,\infty)$, its sample paths have finite variations on $(0,t]$. 
\end{lem}

\subsection{Increment functions of L\'evy processes}

\noindent
For $t \geq 0$, the image of $\omega \in \Omega$ by the rv $Z(t)$ is denoted by $Z(t,\omega)$. For $\omega \in \Omega$, denote by $Z(\cdot,\omega)$ the sample path $t \in [0,\infty) \mapsto Z(t,\omega)$. For $\omega \in \Omega$, the standardized increment function $\eta_{x,\lambda}(\cdot, \omega)$ is defined on $[0,1]$ by 
\begin{equation*}
\eta_{x,\lambda}(s, \omega) = \lambda^{-1}(Z(x+ s\lambda, \omega)-Z(x, \omega)), \enskip s \in [0,1], \quad \textrm{where } x \geq 0, ~\lambda >0. 
\end{equation*}

\noindent
In the sequel, we assume that \textit{all paths} of $Z$ lie in $D([0,\infty))$, which is possible by classical results. Thus, all increment functions of $Z$ lie in $D(0,1)$. The next Lemmas provide measurability and invariance properties.  

\begin{lem}\label{LemMeasEta}
$(i)$ Let $\cB_{\cS}$ be the Borel $\sigma$-algebra of $(D(0,1),\mathcal{S})$. Then, $\cB_{\cS}$ is generated by the class of sets $\mathfrak{C}_{[0,1]}^{D} := D(0,1) \cap \mathfrak{C}_{[0,1]}$.

\noindent
$(ii)$ Let $x \geq 0$ and $\lambda >0$. Then, for any $B \in \cB_{\cS}$, $\left\{ \eta_{x,\lambda} \in B \right\} \in \cF$.

\noindent
$(iii)$ Let $x \geq 0$ and $\lambda >0$. Then, for any $B \in \cB_{\cS}$, we have that 

\noindent
$P\left(\eta_{x,\lambda} \in B \right) = P\left(\eta_{0,\lambda} \in B \right) = P\left( \lambda^{-1}(Z(\cdot ~\lambda)) \in B \right)$.
\end{lem}

\begin{proof}
For $(i)$, see Theorem 12.5 in \cite{billingsley_convergence_1999}. Then, $(ii)$ follows readily from $(i)$. Finally, by Proposition 10.7. in \cite{sato_levy_2013}, 
\begin{equation}\label{invaLevy}
\forall ~B \in \mathfrak{C}_{[0,\infty)}, 
~P\left( Z(x+ \cdot ~\lambda)-Z(x) \in B \right) = 
P\left( Z( \cdot ~\lambda) \in B \right).
\end{equation}

\noindent
For $(iii)$, $(\ref{invaLevy})$ holds for all $B \in \mathfrak{C}_{[0,1]}^{D} \subset \mathfrak{C}_{[0,\infty)}$ and so for all $B \in \sigma\left( \mathfrak{C}_{[0,1]}^{D} \right)=\cB_{\cS}$. 
\end{proof}

\begin{lem}\label{LemMeasSup}
$(i)$ Let $x \geq 0$. For any $u>0$, $F_{x,u} := \left\{ \sup\limits_{0 \leq \tau \leq 1} \left| Z(x+\tau)-Z(x) \right| > u \right\} \in \mathcal{F}$. 

\noindent
$(ii)$ Let $x \geq 0$. For any $u>0$, $P \left( F_{x,u} \right) = P \left( F_{0,u} \right)$.
\end{lem}

\begin{proof}
For all $\omega \in \Omega$, the function $Z(x+ \cdot~, \omega) - Z(x, \omega)$ is right-continuous. So, $F_{x,u} = \left\{ \left( Z(x+ \cdot ~)-Z(x) \right) \in \beta_{u} \right\}$, where $\beta_{u} := \bigcup_{k \in \mathbb{N}} \left\{ f \in D(0,1) : \left|f(\tau_{k})\right| > u \right\}$ and $\left\{\tau_{k} : k \in \mathbb{N}\right\} := [0,1] \cap \mathbb{Q}$. Since $\beta_{u} \in \cB_{\cS}$, we deduce as for $(ii)$ of Lemma $\ref{LemMeasEta}$ that $F_{x,u} \in \mathcal{F}$. This proves $(i)$. Then, by this description of $F_{x,u}$, $(ii)$ follows from $(iii)$ of Lemma $\ref{LemMeasEta}$. 
\end{proof}

\section{Functional spaces}\label{sec:3}

\noindent
In this section, we present properties of the functional spaces considered in the Introduction, and large deviation principles in these spaces. Given a metric space $(\cE,d)$, the open ball of center $f \in \cE$ and radius $\epsilon>0$ is denoted by $B_{d}^{\cE}(f, \epsilon)$.

\subsection{Uniform topology on $D(0,1)$}\label{BasicPropertiesD}

\noindent
On $D(0,1)$, we present the relationships between the Skorohod and the uniform topology $\cU$, induced by the uniform distance denoted by $d_{\cU}$. We denote by $\cU_{0}$ be the $\sigma$-algebra generated by the open balls of $(D(0,1), \cU)$.    

\begin{lem}\label{compareUandS}
$(i)$ The restriction of $\cS$ to $C(0,1)$ coincides there with $\cU$. $(ii)$ $\cU_{0} = \cB_{\cS}$.    
\end{lem}

\begin{proof}
For $(i)$, see Section 12 in \cite{billingsley_convergence_1999}, and for $(ii)$ see Section 15 therein. 
\end{proof}

\vspace{.1cm}
\noindent
Lemma \ref{skoUnif} hereafter is known but a reference for its proof is not available. So, we prove it in subsection \ref{proofLemmaD}.

\begin{lem}\label{skoUnif}
Let $K$ be a compact subset of $(C(0,1),\cU)$ and $\epsilon>0$. Then, there exists $\zeta >0$ such that for all $g \in K$, $B_{\cU}^{D}(g,\epsilon) \supseteq B_{\cS}^{D}(g,\zeta)$, where the subscripts $\cU$ and $\cS$ refer to $d_{\cU}$ and $d_{\cS}$.  
\end{lem}

\subsection{The space $\BV$}\label{BVproperties}

\begin{defi}
For $f \in \BV$, let $|f|_{v}$ be its total variation. For $M>0$, set 
\begin{equation*}
\BVM := \left\{ f \in \BV : |f|_{v} \leq M \right\}.  
\end{equation*}

\noindent
Then, by Proposition 1.4. in \cite{barrio_lectures_2007}, $\BVM$ is a compact subset of $(\BV, \cW)$.
\end{defi}

\noindent
Consider the map $\Lambda : \BV \lra C(0,1)^{*}$ defined by $\langle \Lambda(F) , \phi \rangle := \int_{0}^{1} \phi(x) \mathrm{d}F(x)$, for all $F \in \BV$ and $\phi \in C(0,1)$. By the Riesz representation theorem, $\Lambda$ is a bijection. Then, $C(0,1)^{*}$ is endowed with the weak-$^{*}$ topology and $\BV$ with the topology for which $\Lambda$ is a homeomorphism, called the \textit{weak topology} $\cW$ on $\BV$. Now, $(\BV, \mathcal{W})$ is not metrizable (see Remark $1.2.$ in \cite{barrio_lectures_2007}). However, by results of \cite{hognas1977characterization}, for any $M>0$, the distance $d_{\cH}$ defined hereafter \textit{metricizes} $(\BVM, \cW)$. 
\begin{equation}\label{aa1}
d_{\cH}(f,g)=\int_{0}^{1} |f(u)-g(u)| du + |f(1)-g(1)|, \quad f, g \in \BV. 
\end{equation}

\noindent
Let $\cH$ be the topology on $\BV$ induced by $d_{\cH}$. The $\sigma$-algebra $\cH_{0}$ generated by the open balls of $(\BV, \cH)$, the Borel $\sigma$-algebra $\cB_{\cW}$ of $(\BV, \cW)$ and $\cB_{\cS}$ are linked as follows, which is proved in subsection \ref{proofLemmaBV}.

\begin{lem}\label{H0inclBS}
$(i)$ For $g \in \BV$ and $\xi >0$, the ball $B_{\cH}^{BV} \left( g, \xi \right) := \left\{ f \in BV_{0}(0,1) : d_{\mathcal{H}}(f,g) < \xi \right\}$ is open in $\left( \BV,\cW \right)$.
$(ii)$ As $\sigma$-algebras on $\BV$, we have that $\cH_{0} \subset \cB_{\cW} \subset \cB_{\cS}$. 
\end{lem}

\subsection{Functional large deviations}\label{subsecLDP}

\begin{defi}\label{definitionLDP}
$(i)$ Let $\cE$ be a topological space, with Borel $\sigma$-algebra $\cB$. A function $\cI $ from $\cE$ to $[0, \infty]$ is a rate function if $\cI$ is lower semicontinuous. Furthermore, we say that $\cI$ is a good rate function if its sublevel sets are compact.

\vspace{.1cm}
\noindent
$(ii)$ A family of probability measures $(P_{\lambda})_{\lambda>0}$ on $\left(\cE, \cB \right)$ satisfies a large deviations principle $[LDP]$ in $\cE$, with rate function $\cI$ when, for any closed (resp. open) subset $F$ (resp. $G$) of $\cE$, 
\begin{equation*}
{\underset{\lambda \rightarrow \infty}\varlimsup} \enskip \frac{1}{\lambda} \log P_{\lambda}(F)  \leq - \cI(F) \quad \textrm{and} \quad
{\underset{\lambda \rightarrow \infty}\varliminf} \enskip \frac{1}{\lambda} \log P_{\lambda}(G)  \geq - \mathcal{I}(G) , 
\end{equation*}

\noindent
where for $A \subset \cE, ~\cI(A):=\inf\limits_{f \in A}\cI(f)$. These inequalities are called respectively the upper bound and the lower bound of the LDP.  
\end{defi}

\noindent
Our FERL's under $\left( \cC \right)$ and $\left( \cA \right)$ are based on LDP's for the distributions $(P_{\lambda})_{\lambda>0}$, stated in Theorem $\ref{TheoLDP1}$ and Theorem $\ref{TheoLDP2}$ below, obtained respectively in \cite{varadhan_asymptotic_1966} and \cite{lynch_large_1987}. The rate functions of these LDP's involve the Legendre transform $\Psi$ of the mgf $\Phi$ of $Z(1)$, i.e. for $\alpha \in \R$, $\Psi(\alpha) := \sup \left\{ \alpha \theta - \log \Phi(\theta) : \Phi (\theta) < \infty \right\}$. 

Denote by $AC(0,1)$ the space of absolutely continuous functions on $[0,1]$. Define the map $I : D(0,1) \lra [0, \infty]$ as follows: If $f \in AC(0,1)$ and $f(0)=0$, then $I(f) = \int_{0}^{1} \Psi \left(\frac{\rm df}{{\rm d}s} (s) \right)ds$. Otherwise, $I(f) = \infty$. 

\begin{theo}\label{TheoLDP1}
Under $(\mathcal{C})$, the distributions $(P_{\lambda})_{\lambda > 0}$ satisfy a LDP in $(D(0,1), \cS)$ with good rate function $I$. 
\end{theo}

\noindent
For $f \in BV_0(0,1)$, write $f = f_{+} - f_{-}$, where ${\rm d}f={\rm d}f_{+}-{\rm d}f_-$ is the Hahn-Jordan decomposition of ${\rm d}f$. For $g \in BV_0(0,1)$, write $g = g^{A} + g^{S}$, where ${\rm d}g= {\rm d}g^{A} + {\rm d}g^{S}$ is the Lebesgue decomposition of ${\rm d}g$ into its absolutely continuous and singular components. Then, define the function $J : \BV \lra [0, \infty]$ by
\begin{equation}\label{rateBV}
J(f) = \int_{0}^{1} \Psi \left(\frac{\rm df^{A}}{{\rm d}s} (s) \right) ds + \theta_{0} f_{+}^{S}(1) - \theta_{1} f_{-}^{S}(1), \quad f \in \BV,
\end{equation}  

\noindent
where $\theta_{0} := \sup \left\{ \theta : \Phi(\theta) < \infty \right\} > 0$   and $\theta_{1} := \inf \left\{ \theta : \Phi(\theta) < \infty \right\} < 0$. 

\begin{theo}\label{TheoLDP2}
Let $Z$ be a L\'evy process with generating triplet $\left(A, \nu, \gamma \right)$. If $A=0$ and $\int_{|x| \leq 1} |x|\nu(dx) < \infty$ then, under $(\mathcal{A})$, the distributions $(P_{\lambda})_{\lambda > 0}$ satisfy a LDP in $(\BV, \mathcal{W})$ with good rate function $J$.  
\end{theo}

\section{Main results}\label{sec:4}

\noindent
Due to measurability issues, we consider the class $\cN$ of negligible sets in $\Omega$, that is $\cN := \left\{ N \subset \Omega : \exists N' \in \cF \textrm{ s.t. } N \subset N',~ P(N')=0 \right\}$. Set $\cN^{\rc} := \left\{ N^{\rc} : N \in \cN \right\}$, where $N^{\rc}$ is the complement of $N$ in $\Omega$.

\subsection{Preliminaries}\label{PreliminariesMain}

\noindent
From now, as in the Introduction, $\cE$ is a linear subspace of the space of functions on $[0,1]$ such that, a.s. for all $x \geq 0$ and $\lambda>0$, $\eta_{x,\lambda} \in \cE$. Let $d$ be a distance on $\cE$. For any subset $K$ of $\cE$, consider its $\epsilon$-neighborhood wrt $d$, defined by  
\begin{equation}\label{epaissiBalls}
K^{\epsilon,d} := \left\{ f \in \cE : \exists g \in K \textrm{ s.t. } ~d(f, g)<\epsilon \right\} = \bigcup_{g \in K} B_{d}^{\cE}(g, \epsilon). 
\end{equation}

\noindent
Let $(E_{n})_ {n>1}$ be subsets of $\Omega$. Set $\left\{ E_{n} ~\textrm{ ult.}\right\} := \bigcup\limits_{n>1} \bigcap\limits_{k \geq n} E_{n}$ and $\left\{ E_{n} ~\textrm{ i.o.}\right\} := \bigcap\limits_{n>1} \bigcup\limits_{k \geq n} E_{n}$, where ult. (resp. i.o.) means ultimately (resp. infinitely often).

\begin{defi}\label{FERLdefi}
The Hausdorff distance $\cD$ induced by $d$ is defined as follows. For any subsets $A$ and $B$ of $\cE$, $\cD(A,B) := \inf \left\{ \epsilon>0 : A \subset B^{\epsilon,d} \textrm{ and } B \subset A^{\epsilon,d} \right\}$. Then, the convergence wrt $\cD$ of the sequence $\left( \Ln \right)$ to a subset $\cK \subset \cE$ is equivalent to:
\begin{equation}\label{LnGsubsetK} 
\forall ~\epsilon >0, \quad \left\{ \cL_{n} \subset \cK^{\epsilon;d} ~\textrm{ ult.} \right\} \in \cN^{c}
\quad \textrm{and} \quad 
\left\{ \cK \subset \left(\cL_{n} \right)^{\epsilon;d} ~\textrm{ ult.} \right\} \in \cN^{c}. 
\end{equation}

\noindent
So, $(\ref{LnGsubsetK})$ means that a $FERL$ holds for $(\Ln)$ in $(\cE, d)$ with limit set $\cK \subset \cE$. The first and second parts of $(\ref{LnGsubsetK})$ are called respectively the upper and lower bound in the FERL for $(\cL_{n})$. Indeed, their proofs rely on the upper and the lower bound in some LDP's in $(\cE, \cT)$.
\end{defi}

\noindent
Then, the following result is standard: See for example \cite{barrio_lectures_2007}.

\begin{prop}\label{supFunctional}
If a FERL holds for $\left(\cL_{n} \right)$ in $(\cE,d)$, with limit set $\cK$, then for all continuous map $\Theta : (\cE,d) \longrightarrow \R$, 
\begin{equation}
\left\{ \omega \in \Omega : \sup \left\{\Theta(f) : f \in \cL_{n} \right\}~ {\underset{n \rightarrow \infty} 
\longrightarrow}~ \sup \left\{ \Theta(f) : f \in \cK \right\} \right\} \in \cN^{c}. 
\end{equation}
\end{prop}

\subsection{Results in the general framework}\label{subsecGenFrame}

\noindent
First, Lemma \ref{FromLnNtoLn} provides a general strategy to deduce the FERL for $\left( \Ln \right)$ from that for $\left( \LnN \right)$. Then, in Theorem \ref{generalTheoN}, we give conditions  under which the FERL for $\left( \LnN \right)$ in $(\cE, d)$ holds. Let $\cT$ be a topology on $\cE$ which is not necessarily metricized by the distance $d$. Let $\cB$ be the Borel 
$\sigma$-algebra of $(\cE, \cT)$. Consider the following assumptions:

\vspace{.2cm}
\noindent
$(A1)$: $\cD_{0} \subset \cB \subset \cB_{\cS}$, where $\cD_{0}$ is the $\sigma$-algebra generated by the open balls of $(\cE, d)$.

\vspace{.1cm}
\noindent
$(A2)$: A LDP holds for $(P_{\lambda})_{\lambda>0}$ in $(\cE, \cT)$, with a good rate function denoted by $\cI$.

\vspace{.1cm}
\noindent
$(A3)$: For all $\alpha >0$, there exists a convex subset $E_{\alpha}$ of $\cE$, with $0_{\cE} \in E_{\alpha}$, such that the set $\cK_{\alpha}:=\left\{ f \in \cE : \cI(f) \leq \alpha \right\}$ is included in $E_{\alpha}$ and the distance $d$ metricizes the restriction of $\cT$ to $E_{\alpha}$. 

\vspace{.1cm}
\noindent
$(A4)$: There exists a distance $\delta$ on $\cE$ such that: $(i)$ The topology induced by $\delta$ on $\cE$ is weaker than $\cT$. $(ii)$ For all $\alpha>0$, all compact $K$ of $(E_{\alpha},d)$, $g \in K$ and $\epsilon>0$, there exists $\zeta>0$ such that $B_{d}^{\cE}(g,\epsilon) \supseteq B_{\delta}^{\cE}(g,\zeta)$.

\vspace{.1cm}
\noindent
$(A5)$: For all $\alpha>0$, the function $\cI$ is convex on $E_{\alpha}$ and $\cI(0_{\cE}) = 0$, where $0_{\cE}$ is the null function on $[0,1]$.

\begin{lem}\label{FromLnNtoLn}
Assume that a FERL holds for $\left( \LnN \right)$ in $(\cE, d)$ and condition $(\Delta^{d})$ defined in the Introduction is fulfilled. Then, the FERL holds for $\left( \Ln \right)$ in $(\cE, d)$ with the same limit set as for the FERL for $\left( \LnN \right)$.
\end{lem}

\begin{proof}
The proof is deferred to subsection \ref{ProofLemmaExtend}.
\end{proof}

\begin{theo}\label{generalTheoN}
$(i)$ Assume that $(A1)$-$(A5)$ are fulfilled. Then, the FERL holds for $\left( \LnN \right)$ in $(\cE, d)$ with limit set $\cK_{1/c}$. 

\noindent
$(ii)$ If $(\Delta^{d})$ and $(A1)$-$(A5)$ hold, then the FERL holds for $\left( \Ln \right)$ in $(\cE, d)$, with limit set $\cK_{1/c}$.
\end{theo}

\begin{proof}
First, $(i)$ is proved in subsection \ref{proofGeneralTheo}. Then, $(ii)$ follows from $(i)$ and Lemma \ref{FromLnNtoLn}. 
\end{proof}

\begin{rem}\label{remAssumptions}
Conditions $(A3)$ and $(A4)$ are not as restrictive as they seem. Indeed: 

\vspace{.05cm}
\noindent
$R1)$ If the topology of $d$ coincides with $\cT$, then $(A3)$ and $(A4)$ hold with $E_{\alpha}=\cE$ for all $\alpha>0$, $\delta=d$ and $\zeta=\epsilon$. 

\vspace{.05cm}
\noindent
$R2)$ If the topology of $d$ is weaker than $\cT$, then $(A4)$ holds with $\delta=d$ and $\zeta=\epsilon$. 

\vspace{.05cm}
\noindent
$R3)$ By $(i)$ of Lemma \ref{H0inclBS}, on $\BV$, the topology of $d_{\cH}$ is weaker than $\cW$. So, by $R2)$, $(A4)$ holds in that case. 

\vspace{.05cm}
\noindent
$R4)$ By Lemma \ref{skoUnif}, $(A4)$ holds for $(\cE, \cT, d, \delta)=(D(0,1), \cS, d_{\cU}, d_{\cS})$. 
\end{rem}

\subsection{FERL's under $(\cC)$ and $(\cA)$}

\noindent
For $\alpha>0$, set $K_{\alpha} := \left\{ f \in D(0,1) : I(f) \leq \alpha \right\}$ and $L_{\alpha}:= \left\{ f \in \BV : J(f) \leq \alpha \right\}$, where $I$ and $J$ are the rate functions of subsection \ref{subsecLDP}. We denote by $\id$ the identity function on $[0,1]$ and we set $\mu:=\E[Z(1)]$. 

\begin{theo}\label{ferlPartialCA}
$(i)$ Under $(\cC)$, a FERL holds for $\left( \LnN \right)$ in $(D(0,1), d_{\cU})$, with limit set $K_{1/c}+\mu \id$.

\vspace{.1cm}
\noindent
$(ii)$ Assume that for the generating triplet $\left(A, \nu, \gamma \right)$ of $Z$, $A=0$ and $\int_{|x| \leq 1} |x|\nu(dx) < \infty$. Then, under $(\cA)$, a FERL holds for 
$\left(\LnN \right)$ in $\left(\BV, d_{\cH}\right)$, with limit set $L_{1/c}+\mu \id$. 
\end{theo}

\begin{proof}
See subsection \ref{proofFERLpartial} below. 
\end{proof}

\begin{theo}\label{ferlFullC}
Under $(\cC)$, the FERL holds for $\left( \Ln \right)$ in $(D(0,1), d_{\cU})$, with limit set $K_{1/c}+\mu \id$.
\end{theo}

\begin{proof}
By $(ii)$ of Theorem \ref{generalTheoN}, it is enough to prove that, under $(\cC)$, $(\Delta)$ holds with $d=d_{\cU}$: see subsection \ref{proofFullFERL}.
\end{proof}

\vspace{.2cm}
\noindent
Then, the following Corollary is a consequence of Proposition $\ref{supFunctional}$.

\begin{corr}\label{corrMainTheo}
$(i)$ Under $(\cC)$, for any continuous map $\Theta : (D(0,1),\cU) \lra \R$, 
\begin{equation*}
\sup\left\{ \Theta(\eta_{x,A_{n}}) : x \in [0,n-A_{n}] \right\} ~{\underset{n \rightarrow \infty} \lra }
~\sup\left\{ \Theta(f) : f \in K_{1/c} \right\} \quad \textrm{a.s.}
\end{equation*}

\noindent
Since $\Theta : f \mapsto f(1)$ is such a map, we recover directly the convergence in $(\ref{ERofFrolov})$. 

\vspace{.15cm}
\noindent
$(ii)$ Under the assumptions of $(ii)$ of Theorem $\ref{ferlPartialCA}$, for any continuous map 

\noindent
$\Theta : (\BV, d_{\cH}) \lra \R$, 
\begin{equation*}
\sup\left\{ \Theta(\eta_{x,A_{n}}) : x \in [0,n-A_{n}]  \cap \N \right\} ~{\underset{n \ra \infty} \lra }
~\sup\left\{ \Theta(f) : f \in K_{1/c} \right\} \quad \textrm{a.s.}
\end{equation*}
\end{corr}

\subsection{Proof of Theorem \ref{ferlPartialCA}}\label{proofFERLpartial}

\noindent
We establish Theorem \ref{ferlPartialCA} by applying $(i)$ of Theorem \ref{generalTheoN}. We precise in Table \ref{tbl1} the elements needed to define the framework within which we apply Theorem \ref{generalTheoN}, respectively under $(\cC)$ and under $(\cA)$. 

Namely, under $(\cC)$, we consider $\cE=D(0,1)$, within which all increment functions lie. Under $(\cA)$, we consider $\cE=\BV$, Indeed, by Lemma \ref{pathsBV}, under the assumptions of $(ii)$ of Theorem \ref{ferlPartialCA}, a.s. for all $x \geq 0$ and $\lambda>0$, $\eta_{x,\lambda} \in \BV$. Then, in Table \ref{tbl2}, we indicate the results which prove that Assumptions $(A1)$-$(A4)$ hold for the elements of Table $\ref{tbl1}$. Now, we prove in subsection \ref{reductionCentered} that we may assume that $Z$ is centered, i.e. for all $t \geq 0$, $\E[Z(t)]=0$. In that case, $\mu:=\E[Z(1)]=0$, so that, by Remark \ref{CheckA5}, $(A5)$ holds under $(\cC)$ and $(\cA)$. 

Thus, applying Theorem \ref{generalTheoN}, we obtain that the FERL's of Theorem \ref{ferlPartialCA} hold, with limit set $K_{1/c}$ and $L_{1/c}$. When $Z$ is not centered, by subsection \ref{reductionCentered}, these limit sets are $K_{1/c}+\mu \id$ and $L_{1/c}+\mu \id$. This proves Theorem \ref{ferlPartialCA}.

\begin{table}[h]
\caption{}\label{tbl1}
\begin{tabular}{@{}llllll@{}}
\toprule
$\Phi$ & $\cE$ & $d$ & $\cT$ & $\cI$ & $E_{\alpha}$ \\
\midrule
$(\cC )$ & $D(0,1)$ & $d_{\cU}$ & $\cS$ & $I$ & $AC(0,1)$ \\
\midrule
$(\cA )$ & $BV_{0}(0,1)$ & $d_{\cH}$ & $\cW$ & $J$ & $BV_{0,M_{\alpha}}(0,1)$ \\
\botrule
\end{tabular}
\end{table}

\begin{table}[h]
\caption{}\label{tbl2}
\begin{tabular}{@{}lllll@{}}
\toprule
$\Phi$ & $(A1)$ & $(A2)$ & (A3) & $(A4)$ \\
\midrule
$\left(\cC \right)$ & Lemma \ref{compareUandS}, $(ii)$ & Theorem \ref{TheoLDP1} & Remark \ref{CheckA3}, $(i)$ & Remark \ref{remAssumptions}, $R4)$ \\
\midrule
$\left(\cA \right)$ & Lemma \ref{H0inclBS} & Theorem \ref{TheoLDP2} & Remark \ref{CheckA3}, $(ii)$ & Remark \ref{remAssumptions}, $R3)$ \\
\botrule
\end{tabular}
\end{table}

\begin{rem}\label{CheckA3}
We check that $(A3)$ holds for the elements associated to $(\cC)$ and $(\cA)$ in Table $\ref{tbl1}$.  

\vspace{.1cm}
\noindent
$(i)$ By definition of $I$, for all $\alpha>0$, $K_{\alpha} \subset AC(0,1)$ and by Lemma \ref{compareUandS}, $d_{\cU}$ metricizes the restriction of $\cS$ to $AC(0,1)$. 

\vspace{.1cm}
\noindent
$(ii)$ By (3.9) in \cite{deheuvels1991functional}, for all $\alpha>0$, there exists $M_{\alpha}< \infty$ such that $L_{\alpha} \subset BV_{0, M_{\alpha}}(0,1)$. Then, we mentioned in subsection \ref{BVproperties} that $d_{\cH}$ metricizes the restriction of $\cW$ to $BV_{0, M_{\alpha}}(0,1)$. 
\end{rem}

\begin{rem}\label{CheckA5}
Since $\Psi$ is convex, $I$ is convex on $AC(0,1)$ and $J$ on $\BV$. Moreover, when $\E[Z(1)]=0$, if $f$ is the null function on $[0,1]$, then $I(f)=J(f)=0$. So, $(A5)$ holds for the elements associated to $(\cC)$ and $(\cA)$ in Table $\ref{tbl1}$. 
\end{rem}

\subsection{Examples}

\subsubsection{Continuous paths}

\noindent
Let $\left\{ Z (t) : t \geq 0\right\}$ be a L\'evy process with continuous paths, that is a brownian motion with drift. So, Theorem \ref{ferlFullC} yields the FERL for $\left\{ Z (t) : t \geq 0\right\}$, since it satisfies $(\cC)$.

\subsubsection{Compound Poisson process}

\noindent
Let $\left\{ Y_{i} : i \geq 1 \right\}$ be a sequence of i.i.d. random variables. Let $\left\{ N(t) : t \geq 0 \right\}$ be a homogeneous, right-continuous Poisson process of parameter 
$\lambda$, which is assumed to be independent of $\left\{ Y_{i} : i \geq 1 \right\}$. For any $t \geq 0$, set 
\begin{equation}
S_{N}(t)= \sum\limits_{1 \leq i \leq N(t)} Y_{i}. 
\end{equation}

\noindent
Then, the compound Poisson process $\left\{S_{N}(t) : t \geq 0 \right\}$ is a L\'evy process whose genera\-ting triplet satisfies that $\int_{|x| \leq 1} |x|\nu(dx) \leq \nu(\R) = \lambda < \infty$ and $A=0$. The mgf $\Phi$ is such that for $\theta \in \R$, $\Phi(\theta)=\exp [\lambda( \Phi_{Y_{1}}(\theta)-1)]$, where $\Phi_{Y_{1}}$ is the mgf of $Y_{1}$. So, if 
$\Phi_{Y_{1}}(\theta) < \infty$ for all $\theta$ in $\R$ (resp. in a neighborhood of $0$) then $(\cC)$ (resp. $(\cA)$) holds.

\section{Proofs}\label{sec:5}

\subsection{Proof of Lemma \ref{FromLnNtoLn}}\label{ProofLemmaExtend}

\noindent
Clearly, the lower bound in the FERL for $\left(\LnN \right)$ implies it for $\left(\Ln \right)$. So, Lemma \ref{FromLnNtoLn} follows from Lemma $\ref{generalStrategy}$ hereafter.

\begin{lem}\label{generalStrategy}
Let $\cK$ be a subset of $\cE$. Suppose that, for all $\epsilon>0$ : $(i)$ $\left\{ \LnN \not\subset \cK^{\epsilon;d} ~\mathrm{ i.o.} \right\} \in \cN$ and 
$(ii)$ $\left\{ \Delta_{n}^{d}(\epsilon) ~\mathrm{ i.o.} \right\} \in \cN$. Then, for all $\epsilon>0$, $\left\{ \Ln \not\subset \cK^{\epsilon;d} ~\mathrm{ i.o.} \right\} \in \cN$.
\end{lem}

\begin{proof}
Let $\epsilon>0$. For $n>1$, set $\Lambda_{n}^{d}(\epsilon) :=
\left\{ \Ln \not\subset \cK^{\epsilon;d} \right\} \bigcap \left\{ \LnN \subset \cK^{\epsilon/2;d} \right\}$. Thus,   
\begin{align}
\Lambda_{n}^{d}(\epsilon)
&\subset
\bigcup_{x \in [0, n-A_{n}] \cap \N^{c}} \left[ \left\{\eta_{x,A_{n}} \notin \cK^{\epsilon;d} \right\} \bigcap \left\{ \LnN \subset \cK^{\epsilon/2;d} \right\} \right] \\
\label{2ndLineLambda} & \subset
\bigcup_{x \in [0, n-A_{n}] \cap \N^{c}} \left[ \left\{\eta_{x,A_{n}} \notin \cK^{\epsilon;d} \right\} \bigcap \left\{ \eta_{\lfloor x \rfloor,A_{n}} \in \cK^{\epsilon/2;d} \right\} \right] \\
\label{3rdLineLambda} & \subset
\left\{ \exists x \in [0, n-A_{n}] \textrm{ ~s.t.~ } d\left( \eta_{x,A_{n}} , \eta_{\lfloor x \rfloor,A_{n}} \right) \geq \epsilon/2 \right\} =: \Delta_{n}^{d}(\epsilon/2). 
\end{align}

\noindent
Indeed, if $\eta_{\lfloor x \rfloor,A_{n}} \in \cK^{\epsilon/2;d}$ and $\eta_{x,A_{n}} \notin \cK^{\epsilon;d}$ as in $(\ref{2ndLineLambda})$, then for some $g \in \cK$, we have that $d\left( \eta_{\lfloor x \rfloor,A_{n}} , g \right) < \frac{\epsilon}{2}$ and $d\left( \eta_{x,A_{n}} , g \right) \geq \epsilon$. Therefore, by the triangle inequality, $d\left(\eta_{x,A_{n}} , \eta_{\lfloor x \rfloor,A_{n}}\right) \geq \frac{\epsilon}{2}$. This poves that $(\ref{2ndLineLambda})$ implies $(\ref{3rdLineLambda})$. Now, 
\begin{equation}\label{LnLambda}
\left\{ \Ln \not\subset \cK^{\epsilon;d} ~\textrm{ i.o.} \right\} \bigcap 
\left\{ \LnN \subset \cK^{\epsilon/2;d} ~\textrm{ ult.} \right\} \subset
\left\{ \Lambda_{n}^{d}(\epsilon) \textrm{ i.o.} \right\}.
\end{equation}

\noindent
Recall from $(\ref{3rdLineLambda})$ that $\left\{ \Lambda_{n}^{d}(\epsilon) \textrm{ i.o.} \right\} \subset \left\{ \Delta_{n}^{d}(\epsilon/2) \textrm{ i.o.} \right\}$. So, 
 $(\ref{LnLambda})$ then $(ii)$ and $(i)$ imply that 
\begin{equation*}
\left\{ \Ln \not\subset \cK^{\epsilon;d} ~\textrm{i.o.} \right\} \subset \left\{ \Lambda_{n}^{d}(\epsilon) \textrm{ i.o.} \right\} \bigcup 
\left\{ \LnN \not\subset \cK^{\epsilon/2;d} ~\textrm{ i.o.} \right\} \in \cN.  
\end{equation*}
\end{proof}

\subsection{Proof of Theorem \ref{generalTheoN}}\label{proofGeneralTheo}

\noindent
The following Remark allows to reindex series to apply the Borel-Cantelli lemma. 

\begin{rem}\label{njReduce}
For any integer $j >1$, set $n_{j} := \max \left\{n : A_{n}=j \right\}$, so that 
\begin{equation}\label{njBounds}
\exp \left( j/c \right) \leq n_{j}<\exp \left( (j+1)/c \right). 
\end{equation}

\noindent
If $n_{j} < n \leq n_{j+1}$, then $A_{n}=A_{n_{j+1}}=j+1$. So, for such $n$, $n-A_{n} \leq n_{j+1}-A_{n_{j+1}}$ and for all $x \in [0,n-A_{n}] \subset [0,n_{j+1}-A_{n_{j+1}}]$,  
$\eta_{x,A_{n}}=\eta_{x,A_{n_{j+1}}}$. In particular, $\Ln \subset \cL_{n_{j+1}}$. 
\end{rem}

\noindent
Then, we prove Theorem \ref{generalTheoN} as follows. We only need to establish $(i)$. \\

\begin{proof}
We prove the upper and then, the lower bound in the FERL for $\left(\LnN \right)$ in $(\cE, d)$, using Lemma \ref{LemLDP1} below as an auxiliary result. First, by $(A3)$ and $(A4)$, for all $\alpha>0$, $\cK_{\alpha}$ is compact in $(E_{\alpha}, d)$, in $(\cE, d)$ and in $(\cE, \delta)$. Fix $\epsilon >0$. 

\vspace{.2cm}

\textit{Upper bound}: First, we treat measurability issues. By $(ii)$ of $(A4)$ and $(\ref{epaissiBalls})$, there exists $\zeta >0$ such that $\left(\cK_{1/c}\right)^{\epsilon; d} \supseteq \left(\cK_{1/c}\right)^{\zeta; \delta}$. Since $\cK_{1/c}$ is separable in $(\cE,d)$ and in $(\cE, \delta)$, Lemma 1 in Section 6 of \cite{billingsley_convergence_1999} implies that $\left(\cK_{1/c}\right)^{\epsilon; d} \in \cD_{0}$ and $\left(\cK_{1/c}\right)^{\zeta; \delta}$ lies in the $\sigma$-algebra generated by the open balls wrt $\delta$, which is included in $\cB$, by $(i)$ of $(A4)$. So, by $(A1)$ and $(ii)$ of Lemma $\ref{LemMeasEta}$, for all $n>1$, $\left\{ \LnN \subset \left(\cK_{1/c}\right)^{\epsilon; d} \right\} \in \cF$ and 
$\left\{ \LnN \subset \left(\cK_{1/c}\right)^{\zeta; \delta} \right\} \in \cF$. Then,  
\begin{equation}\label{firstSteProof}
P \left(\LnN \not\subset \left(\cK_{1/c}\right)^{\epsilon; d} \right) \leq 
P \left(\LnN \not\subset \left(\cK_{1/c}\right)^{\zeta; \delta} \right) = 
\sum\limits_{m=0}^{n-A_{n}} P\left( \eta_{m, A_{n}} \in \left[ \left(\cK_{1/c}\right)^{\zeta; \delta} \right]^{\rc} \right). 
\end{equation}

\noindent
So, by $(\ref{firstSteProof})$ and $(iii)$ of Lemma \ref{LemMeasEta}, since $\eta_{0, A_{n}}=A_{n}^{-1}Z(\cdot ~A_{n})$,
\begin{equation*}
P \left(\LnN \not\subset \left(\cK_{1/c}\right)^{\epsilon; d} \right) \leq (n-A_{n}+1) P\left( A_{n}^{-1}Z(\cdot ~A_{n}) \in \left[ \left(\cK_{1/c}\right)^{\zeta; \delta} \right]^{\rc} \right). 
\end{equation*}

\noindent
Now, $\cK_{1/c}$ is compact in $(E_{1/c}, d)$. So, by $(\ref{epaissiBalls})$ and $(i)$ of $(A4)$, $\left[ \left(\cK_{1/c}\right)^{\zeta; \delta} \right]^{\rc}$ is $\cT$-closed. Then, by $(A2)$, for all $\theta >0$, there exists $N_{\theta} \in \N$ such that for all $n \geq N_{\theta}$,
\begin{equation} \label{m} 
 \enskip P\left( A_{n}^{-1}Z(\cdot ~A_{n}) \in \left[ \left(\cK_{1/c}\right)^{\zeta; \delta} \right]^{c} \right) \leq 
\exp \left[ A_{n} \left( -\cI_{1/c}^{\zeta} + \theta \right) \right], 
\end{equation}

\noindent
where $\cI_{1/c}^{\zeta} := \inf \left\{ \cI(x) : x \notin  \left(\cK_{1/c}\right)^{\zeta; \delta} \right\}$. Now, $(i)$ of Lemma $\ref{LemLDP1}$ implies that $\cI_{1/c}^{\zeta} = 1/c + \alpha$ with $\alpha >0$. So, for all $n \geq N_{\theta}$, 
\begin{equation}\label{n}
P \left( \LnN \not\subset \left(\cK_{1/c}\right)^{\epsilon; d} \right) \leq n \exp \left[ A_{n} \left( - \frac{1}{c} - \frac{3\alpha}{4} \right) \right].
\end{equation}

\noindent
where we applied $(\ref{m})$ with $\theta = \frac{\alpha}{4}$. By $(\ref{n})$ applied with $n=n_{j}$, so that $A_{n}=j$, and $(\ref{njBounds})$, 
\begin{equation}\label{endUpperBound}
P \left(\cL_{n_{j}}^{\N} \not\subset \left(\cK_{1/c}\right)^{\epsilon; d} \right) \leq
n_{j} \exp \left[ -j \left( \frac{1}{c} + \frac{3\alpha}{4} \right) \right] < 
\exp \left( \frac{1}{c}-j \frac{3\alpha}{4}  \right).
\end{equation}

\noindent
So the upper bound in the FERL for $\left( \LnN \right)$ is proved. Indeed, by Remark \ref{njReduce}, $(\ref{endUpperBound})$ and the Borel-Cantelli lemma,  
\begin{equation*}
\left\{ \LnN \not\subset \left(\cK_{1/c}\right)^{\epsilon; d} ~\mathrm{ i.o.} \right\} \subset \left\{ \cL_{n_{j}}^{\N} \not\subset \left(\cK_{1/c}\right)^{\epsilon; d} ~\mathrm{ i.o. ~in ~j} \right\} \in \cN.
\end{equation*}

\vspace{.15cm}

\textit{Lower bound}: We fix $g \in \cK_{1/c}$. Then, for any $n>1$, we set $R_{n}:=[(n-A_{n})/A_{n}]$ and $\cQ_{n} := \left\{ \eta_{rA_{n}, A_{n}} : 0 \leq r \leq R_{n}-1 \right\}$. Now, 
\begin{equation}\label{firstLB}
\left\{ g \notin (Q_{n})^{\epsilon/2;d} \right\} = \bigcap_{r=0}^{R_{n}-1} \left\{\eta_{rA_{n}, A_{n}} \notin B_{d}^{\cE}(g,\epsilon/2) \right\} \subset 
\bigcap_{r=0}^{R_{n}-1} \left\{\eta_{rA_{n}, A_{n}} \notin B_{\delta}^{\cE}(g,\zeta') \right\},  
\end{equation}

\noindent
for $\zeta'>0$, by $(ii)$ of $(A4)$. As for the upper bound, $\left\{\eta_{rA_{n}, A_{n}} \notin B_{d}^{\cE}(g,\epsilon/2) \right\} \in \cF$ and $\left\{\eta_{rA_{n}, A_{n}} \notin B_{\delta}^{\cE}(g,\zeta') \right\} \in \cF$, for all $0 \leq r \leq R_{n}-1$. Thus, $(\ref{firstLB})$, the mutual independence of the $\left(\eta_{rA_{n}, A_{n}} \right)_{r}$ and $(iii)$ of Lemma \ref{LemMeasSup} imply that 
\begin{equation*}
P\left( g \notin (\cQ_{n})^{\epsilon/2;d} \right) \leq
\prod_{r=0}^{R_{n}-1} \left[ 1-P\left( \eta_{rA_{n}, A_{n}} \in B_{\delta}^{\cE}(g,\zeta') \right) \right] =
\left[ 1-P \left( A_{n}^{-1}Z(\cdot ~A_{n}) \in B_{\delta}^{\cE}(g,\zeta') \right)  \right]^{R_{n}}.
\end{equation*}

\noindent
Now, $\cK_{1/c}$ is compact in $(E_{1/c}, d)$. Then, by $(i)$ of $(A4)$, $B_{\delta}^{\cE}(g,\zeta')$ is $\cT$-open. So, by $(A2)$, for all $\theta>0$, there exists $N_{\theta} \in \N$ such that for all $n \geq N_{\theta}$, 
\begin{equation}\label{LDPlowC}
P\left( A_{n}^{-1}Z(\cdot ~A_{n}) \in B_{\delta}^{\cE}(g,\zeta') \right) \geq 
\exp \left( A_{n}\left(-\cI \left( B_{\delta}^{\cE}(g,\zeta') \right)-\theta \right) \right). 
\end{equation}

\noindent
By $(ii)$ of Lemma $\ref{LemLDP1}$, $I\left( B_{\delta}^{\cE}(g,\zeta')  \right) = 1/c - \beta$ with $\beta >0$. Then, we apply $(\ref{LDPlowC})$ with $\theta = \frac{\beta}{4}$. So, 
for all $n \geq N_{\theta}$, 
\begin{equation}\label{BClow}
P \left( g \notin (Q_{n})^{\epsilon/2;d} \right) \leq \left[ 1 - \exp \left( A_{n} \left( -\frac{1}{c} + \frac{3\beta}{4}  \right) \right) \right]^{R_{n}} \textrm{~i.e.~}
 P \left( g \in (Q_{n})^{\epsilon/2; d} ~\textrm{ ult.} \right) = 1.
\end{equation}

\noindent
Since $\cK_{1/c}$ is compact in $\left( \cE, d \right)$, there exist functions $(g_{q})_{q=1,...,d}$ which belong to $\cK_{1/c}$ such that $\cK_{1/c} \subset \bigcup_{q=1}^{d} B_{d}^{\cE}(g_{q},\epsilon /2)$. So,  
\begin{equation*}
\forall ~n>1, \enskip \left\{ \left\{ g_{q} : q=1,...,d \right\} \subset (Q_{n})^{\epsilon/2; d} \right\} \subset \left\{ \cK_{1/c} \subset \left(Q_{n}\right)^{\epsilon;d} \right\}, 
\end{equation*}

\noindent
which holds by the triangle inequality. By $(\ref{BClow})$ applied to each $g_{q}$, $P\left( \left\{ g_{q} : q=1,...,d \right\} \subset (Q_{n})^{\epsilon/2;d} ~\textrm{ ult.}
\right) = 1$. So, 
\begin{equation*}
\left\{ \cK_{1/c} \subset \left(\mathcal{L}_{n}\right)^{\epsilon; d}~\textrm{ ult.} \right\} \supset 
\left\{ \cK_{1/c} \subset \left(Q_{n}\right)^{\epsilon; d} ~\textrm{ ult.} \right\}
\in \mathcal{N}^{c}. 
\end{equation*}

\noindent
So, the lower bound of the FERL for $\left( \LnN \right)$ and thus the FERL for $\left( \LnN \right)$ are proved. 
\end{proof}

\begin{lem}\label{LemLDP1}
$(i)$ For all $\alpha>0$ and $\zeta>0$, $\cI_ {\alpha}^{\zeta} := \inf \left\{ \cI(f) : f \notin (\cK_{\alpha})^{\zeta; \delta} \right\} > \alpha$.

\vspace{.1cm}
\noindent
$(ii)$ For all $\alpha>0$, $g \in \cK_{\alpha}$ and $\zeta>0$, $\cI \left( B_{\delta}^{\cE}(g,\zeta) \right) := \inf \left\{ \cI(f) : f \in B_{\delta}^{\cE}(g, \zeta) \right\} < \alpha$.

\end{lem}

\begin{proof}
For $(i)$, the assumptions imply that $\cI$ is a good rate function on $(\cE, \delta)$. Then, we mimic the proof of Lemma 2.1 in \cite{sanchis_functional_1994} to prove $(i)$. For $(ii)$, let $g \in \cK_{\alpha}$ and $\zeta>0$. The map $\phi : [0,1] \lra [0, \infty)$ defined by $\phi(\lambda)=\delta(\lambda g; g)$ is continuous and $\phi(\lambda)=0$ iff $\lambda=1$. So, there exists $\lambda_{0}<1$ such that $\phi(\lambda_{0})<\zeta$ i.e. $\lambda_{0}g \in B_{\delta}^{\cE}(g, \zeta)$. Now, by $(A3)$, $E_{\alpha}$ is convex so that $\lambda_{0}g \in E_{\alpha}$. Then, by $(A5)$, $\cI(\lambda_{0}g) \leq \lambda_{0} \cI(g) < \cI(g) \leq \alpha$, which proves $(ii)$. 
\end{proof}

\subsection{Reduction to the centered case}\label{reductionCentered}

\noindent
Here, we prove that for Theorem $\ref{ferlPartialCA}$ and Theorem $\ref{ferlFullC}$, one can suppose that $Z$ is centered.

\begin{lem}\label{translation}
For $\mu \in \R$, define the L\'evy process $Z^{\mu}$ by $Z^{\mu}(t) = Z(t) + \mu t$, $t \geq 0$. Any increment function of $Z^{\mu}$ is denoted by $\eta_{x,A_{n}}^{\mu}$. Let $I^{\mu}$ and $J^{\mu}$ be the analogues for $Z^{\mu}$ of the rate functions $I$ and $J$. Set $\left( \cL_{n}^{\mu} \right) := \left\{ \eta_{x,A_{n}}^{\mu}: x \in [0, n-A_{n}] \right\}$ and for any $\alpha>0$, set $K_{\alpha}^{\mu} := \left\{ f \in D(0,1) : I^{\mu}(f) \leq \alpha \right\}$ and $L_{\alpha}^{\mu} := \left\{ f \in \BV : J^{\mu}(f) \leq \alpha \right\}$. 

\vspace{.2cm}

\noindent
Hereunder, the metric space $(\cE, d)$ designates either $(D(0,1), d_{\cU})$ or $(\BV, d_{\cH})$. Then, for all $\epsilon>0$, 
\begin{equation}\label{eqTranslationFinal}
\left( \cL_{n}^{\mu} \right)^{\epsilon,d} = \left( \cL_{n} \right)^{\epsilon,d} + \mu \id 
\quad ; \quad 
\left( K_{\alpha}^{\mu} \right)^{\epsilon,d} = \left( K_{\alpha} \right)^{\epsilon,d} + \mu \id
\quad \textrm{and}\quad 
\left( L_{\alpha}^{\mu} \right)^{\epsilon,d} = \left( L_{\alpha} \right)^{\epsilon,d} + \mu \id.
\end{equation}
\end{lem}

\begin{proof}
The proof follows readily from the definitions. We omit it.  
\end{proof}

\begin{lem}\label{ZtIstZ1}
If $(\cA)$ holds then, for all $t \geq 0$, $\E[Z(t)] = t\E[Z(1)]$. 
\end{lem}

\begin{proof}
Under $(\mathcal{A})$, $\E[|Z(1)|] < \infty$. Thus, by Theorem 25.3 in \cite{sato_levy_2013}, for all $t \geq 0$, $\E[|Z(t)|] < \infty$. Denoting by $\widehat{\mu}$ the characteristic function of $Z(1)$, that of $Z(t)$ is $(\widehat{\mu})^{t}$. We deduce that $\E[Z(t)] = \frac{1}{i} \frac{\partial (\widehat{\mu}(z))^{t}}{\partial z} \big|_{z=0} = t \E[Z(1)]$. 
\end{proof}

\vspace{.1cm}
In conclusion, by Lemma \ref{ZtIstZ1}, given a L\'evy process $Z$ such that $(\cA)$ holds, the L\'evy process $\overline{Z}$ defined by $\overline{Z}(t) = Z(t)-t\E[Z(1)]$, $t \geq 0$, is centered. Then, by Lemma $\ref{translation}$ applied with $\mu = \E[Z(1)]$, if the FERL holds for the increments of $\overline{Z}$ with limit set $K_{1/c}$ or $L_{1/c}$, then it holds for those of $Z$ with limit set $K_{1/c}+\mu \id$ or $L_{1/c}+\mu \id$.

\subsection{Proof of Theorem \ref{ferlFullC}}\label{proofFullFERL}

\noindent
We prove in Lemma \ref{UppMainLem2} that, under $(\cC)$, $(\Delta^{d_{\cU}})$ holds, for which Lemma \ref{UpperPreLem2} hereafter is a preliminary result. Recall from subsection \ref{reductionCentered} that we assume that $Z$ is centered.

\begin{lem}\label{UpperPreLem2}
For all $u>0$ and $n>1$, $\Delta_{n}^{d_{\cU}}(u) \subset \bigcup\limits_{i=0}^{n+1} \left\{ \sup\limits_{0 \leq \tau \leq 1}\left|Z(i+\tau)-Z(i)\right| > 
\frac{uA_{n}}{9} \right\}$. 

\end{lem}

\begin{proof}
By the triangle inequality,  
\begin{align*}
\Delta_{n}^{d_{\cU}}(u) 
&\subset 
\bigcup\limits_{x \in [0, n-A_{n}]}
\left\{ \left\| Z(x+\cdot ~A_{n})-Z(\lfloor x \rfloor+\cdot ~A_{n})\right\|_{\mathcal{U}} + \left\| Z(x)-Z(\lfloor x \rfloor)\right\|_{\cU} \geq uA_{n} \right\} \\
&\subset
\bigcup\limits_{x \in [0, n-A_{n}]}
\left[ \left\{ \left\| Z(x+\cdot ~A_{n})-Z(\lfloor x \rfloor+\cdot ~A_{n})
\right\|_{\mathcal{U}} > \frac{uA_{n}}{3} \right\} \bigcup
\left\{ \left| Z(x)-Z(\lfloor x \rfloor)\right| > \frac{uA_{n}}{3} \right\} \right] \\
&=  
\bigcup\limits_{x \in [0, n-A_{n}]}
\left[ \bigcup\limits_{s \in [0,1]}
\left\{ \left| Z(x+sA_{n})-Z(\lfloor x \rfloor+sA_{n})
\right| > \frac{uA_{n}}{3} \right\} \right] \\
&\subset
\bigcup\limits_{y \in [0,n]} 
\left[ \bigcup\limits_{a \in [0,1]} \left\{ \left|Z(y+a)-Z(y)\right| > \frac{uA_{n}}{3} \right\} \right]. 
\end{align*}

\noindent
Now, for any $y \in [0,n]$ and $a \in [0,1]$, two cases occur.

\noindent\\
\textit{First case}: If $y+a \leq \lfloor y \rfloor+1$, then 
\begin{equation*}
|Z(y+a)-Z(y)| \leq  \left| Z(y+a)-Z\left(\lfloor y \rfloor \right) \right| + \left| Z(y)-Z\left(\lfloor y \rfloor \right) \right| 
\leq
2 \sup\limits_{0 \leq \tau \leq 1} 
\left| Z(\lfloor y \rfloor + \tau)-Z\left(\lfloor y \rfloor \right) \right|.
\end{equation*}

\noindent
\textit{Second case}: If $y+a > \lfloor y \rfloor+1$, then 
$0 < (y+a)-(\lfloor y \rfloor+1) < (y+1)-y = 1$, so that  
\begin{align*}
|Z(y+a)-Z(y)| &\leq
\left| Z(y+a)-Z\left( \lfloor y \rfloor+1 \right) \right| +
\left| Z\left(\lfloor y \rfloor+1\right)-Z\left(\lfloor y \rfloor\right) \right| + \left| Z(y)-Z\left(\lfloor y \rfloor\right) \right| \\
&\leq
\sup\limits_{0 \leq \tau \leq 1}
\left| Z(\lfloor y \rfloor+1 + \tau)-Z\left(\lfloor y \rfloor+1 \right) \right| +
2 \sup\limits_{0 \leq \tau \leq 1} 
\left| Z(\lfloor y \rfloor+ \tau)-Z\left(\lfloor y \rfloor \right) \right|. 
\end{align*}

\noindent
So, we conclude, since we obtain in both cases that  
\begin{equation*} \label{max} 
\sup\limits_{y \in [0,n]} 
\left\{ \sup\limits_{0 \leq a \leq 1} |Z(y+a)-Z(y)| \right\} \leq
3 \max\limits_{0 \leq i \leq n+1} 
\left\{ \sup\limits_{0 \leq \tau \leq 1} \left| Z\left(i + \tau \right)-Z(i) 
 \right| \right\}.
\end{equation*}
\end{proof}

\begin{lem}\label{UppMainLem2}
Assume that $(\cC)$ holds. Then, for all $\epsilon >0$, $\left\{ \Delta_{n}^{d_{\cU}}(\epsilon) ~\mathrm{ i.o.} \right\} \in \cN$. 
\end{lem}

\begin{proof}
By Lemma \ref{LemMeasSup}, for all $n>1$ and $u>0$, $\bigcup\limits_{i=0}^{n+1} \left\{ \sup\limits_{0 \leq \tau \leq 1}\left|Z(i+\tau)-Z(i)\right| > \frac{u A_{n}}{9} \right\} \in \mathcal{F}$. Set    
\begin{equation*}
\pi_{n}(u) := P\left( \bigcup\limits_{i=0}^{n+1} \left\{ \sup\limits_{0 \leq \tau \leq 1}\left|Z(i+\tau)-Z(i)\right| > \frac{u A_{n}}{9} \right\} \right).
\end{equation*}

\noindent
Fix $\epsilon>0$. By $(ii)$ of Lemma \ref{LemMeasSup}, for all $j$ large enough, 
\begin{equation*}
\pi_{n_{j}}(\epsilon) \leq \sum\limits_{i=0}^{n_{j}+1} P\left( \sup\limits_{0 \leq \tau \leq 1} \left|Z(i+\tau)-Z(i)\right| > \frac{j \epsilon}{9} \right) \leq 
2n_{j} P\left( \sup\limits_{0 \leq \tau \leq 1} \left|Z(\tau)\right| > \frac{j \epsilon}{9} \right).
\end{equation*}

\noindent
Now, for all $\theta >0$ and $v>0$,  
\begin{equation*}
P\left( \sup\limits_{0 \leq \tau \leq 1} \left|Z(\tau)\right| > v \right)  
\leq
P\left( \sup\limits_{0 \leq \tau \leq 1} \exp(\theta Z(\tau)) > \exp(\theta v) \right) +
P\left( \sup\limits_{0 \leq \tau \leq 1} \exp(-\theta Z(\tau)) > \exp(\theta v) \right).
\end{equation*}

\noindent
Since $Z$ is centered, $Z$ is a martingale. So, the processes $\left\{ \exp[\theta Z(t)] : t \geq 0 \right\}$ and $\left\{ \exp[-\theta Z(t)] : t \geq 0 \right\}$ are nonnegative submartingales. Then, $(\cC)$ implies that for all $\theta >0$, $\Phi(\theta)+\Phi(-\theta)< \infty$. By Doob's inequality, for all $\theta>0$, $j>1$, 
\begin{equation}\label{piNj}
\pi_{n_{j}}(\epsilon) \leq 2n_{j}\left[\Phi(\theta)+\Phi(-\theta)\right] \exp\left(-\theta j\frac{\epsilon}{9}\right) < 
2\left[\Phi(\theta)+\Phi(-\theta)\right]\exp\left(\frac{1}{c}\right) \exp\left(-j\left(\theta\frac{\epsilon}{9}-\frac{1}{c} \right) \right).
\end{equation}

\noindent
Indeed, by definition, $n_{j} < \exp\left(\frac{j+1}{c}\right)$. By 
$(\mathcal{C})$, we can choose $\theta = \theta(\epsilon,c)$ such that 
\begin{equation}\label{epsilonChoice}
\theta \frac{\epsilon}{9}>\frac{1}{c}
\quad \textrm{and} \quad
\Phi(\theta)+\Phi(-\theta)< \infty. 
\end{equation}

\noindent
By $(\ref{piNj})$ and $(\ref{epsilonChoice})$, $\sum\limits_{j>1} \pi_{n_{j}}(\epsilon) < \infty$. By Lemma \ref{UpperPreLem2} and the Borel-Cantelli lemma,
\begin{equation*}
\left\{\Delta_{n_{j}}^{\cU}(\epsilon) \textrm{ i.o. in } j \right\} \subset \left\{ \bigcup\limits_{i=0}^{n_{j}+1} \left\{ \sup\limits_{0 \leq \tau \leq 1}\left|Z(i+\tau)-Z(i)\right| > \frac{j\epsilon}{9} \right\} \textrm{ i.o. in } j\right\} \in \mathcal{N}. 
\end{equation*}

\noindent
We conclude by Remark $\ref{njReduce}$, from which $\left\{ \Delta_{n}^{d_{\cU}}(\epsilon) ~\mathrm{ i.o.} \right\} \subset \left\{ \Delta_{n_{j}}^{d_{\cU}}(\epsilon) ~\mathrm{ i.o. ~ in ~ j} \right\} \in \cN$.
\end{proof}

\subsection{Proof of Lemma \ref{skoUnif}}\label{proofLemmaD}

\begin{proof}
For any $f \in C(0,1)$, let $\omega_{f}$ be the modulus of continuity of $f$, defined for $\delta >0$ by
\begin{equation*}
\omega_{f} (\delta) = \sup \left\{ \left|f(x)-f(y)\right| : \left|x-y\right| \leq \delta \right\}.
\end{equation*}

\noindent
By the Arzel\`a-Ascoli theorem, for all $\epsilon>0$, there exists 
$\delta_{\epsilon} > 0$ such that
\begin{equation}\label{moduleC}
\sup \left\{ \omega_{g}(\delta_{\epsilon}) : g \in K \right\} < \epsilon/2.
\end{equation}

\noindent
Let $g \in K$ and $\epsilon>0$. Set $\zeta :=\min \left\{ \delta_{\epsilon} ; \epsilon/2 \right\}$. By $(\ref{defsko})$, for all $h \in B_{\cS}^{D}(g,\zeta)$, there exists $\nu_{h} \in \Lambda$ satisfying
\begin{equation}\label{nuH}
\left\| \nu_{h}-I \right\| < \zeta \leq \delta_{\epsilon} 
\quad \textrm{and} \quad  
\left\| h-g\circ \nu_{h} \right\| < \zeta \leq \epsilon/2.
\end{equation}

\noindent
The first part of $(\ref{nuH})$ combined to $(\ref{moduleC})$ imply that $\left\| g\circ \nu_{h}-g \right\| \leq w_{g}(\left\| \nu_{h}-I \right\|) < \epsilon/2$, which combined to the second part of $(\ref{nuH})$ implies that $\left\| h-g \right\| \leq  \left\| h-g\circ\nu_{h} \right\| +\left\| g\circ\nu_{h}-g \right\| < \epsilon$. So, $h \in B_{\cU}^{D}(g,\epsilon)$. 
\end{proof}

\subsection{Proof of Lemma \ref{H0inclBS}}\label{proofLemmaBV}

\noindent
Since $\left( \BV,\cW \right)$ is not metrizable, we characterize its topology by  convergence of nets. By results of \cite{hognas1977characterization}, a net $\left(f_{\alpha}\right)$ converges to $\widehat{f}$ in $\left( \BV,\cW \right)$ iff there exists $M>0$ such that $f_{\alpha}$ is ultimately in $\BVM$ and $d_{\cH} \left( f_{\alpha}, \widehat{f} \right) \rightarrow 0$. Now, we may prove $(i)$ and $(ii)$ of Lemma $\ref{H0inclBS}$. 

\vspace{.2cm}

For $(i)$, let $g \in \BV$ and $\theta_{g} : \left( \BV, \cW \right) \lra \R$ be the map defined by $\theta_{g}(f) = d_{\cH}(f,g)$. Let $(f_{\alpha})$ be a net in $\BV$ weakly converging to $\widehat{f}$. Then, $\left| d_{\cH}\left(f_{\alpha},g \right) - d_{\cH}\left(\widehat{f},g \right) \right| \leq d_{\cH}(f_{\alpha},\widehat{f}) \rightarrow 0$. So, the net $\left( \theta_{g}\left( f_{\alpha} \right) \right)$ converges to $\theta_{g}\left( \widehat{f} \right)$ and thus the map $\theta_{g}$ is continuous. This proves $(i)$. The following Lemma proves $(ii)$.

\begin{lem}\label{SkorohodWeak}
$(i)$ For any $M>0$, on $\BVM$, $\cS$ is stronger than $\cH$. $(ii)$ For any $M>0$, $BV_{0,M}(0,1)$ is a closed subset of $(D(0,1), \cS)$. $(iii)$ We have that $\cB_{\cW} \subset \cB_{\cS}$. 
\end{lem}
 
\begin{proof}
Let $(f_{n})_{n \geq 1}$ be a sequence of functions in $D(0,1)$ and $f \in D(0,1)$. We prove $(i)$ by $(\ref{CVskoCVhog})$ below, written in \cite{deheuvels1991functional}. 
\begin{equation}\label{CVskoCVhog}
d_{\cS}(f_{n},f) {\underset{}\longrightarrow} 0 ~\textrm{ as } n \rightarrow \infty
\enskip \implies \enskip
d_{\cH}(f_{n},f) {\underset{}\longrightarrow} 0 ~\textrm{ as } n \rightarrow \infty.
\end{equation}

\noindent
For $(ii)$, let $(f_{n})$ be a sequence in $\BVM$ converging to $f$ in $(D(0,1), \cS)$. Now, $\BV$ is compact in $(\BV, \cW)$. So, there exists a subsequence $(f_{\phi(n)})$ which converges weakly to some $\widetilde{f} \in BV_{0,M}(0,1)$, so that $d_{\cH}(f_{\phi(n)},\widetilde{f}) {\underset{}\longrightarrow} 0$. Now, by $(\ref{CVskoCVhog})$, $d_{\cH}(f_{\phi(n)},f) {\underset{}\longrightarrow} 0$. So, $\widetilde{f} = f \in BV_{0,M}(0,1)$, which proves $(ii)$. 

For $(iii)$, we prove that for all closed subset $\Gamma$ of $(\BV, \cW)$ and that for all $M \in \N$, $\Gamma_{M}:=\Gamma \cap \BVM \in \mathcal{B}_{\cS}$. Now, for all $M \in \N$, $\Gamma_{M}$ is a closed subset of $(BV_{0,M}(0,1), \cW)$. Since $d_{\cH}$ metricizes $(\BVM, \cW)$, $(i)$ implies that $\Gamma_{M}$ is a closed subset of $(BV_{0,M}(0,1), \cS)$ and so of $(D(0,1), \cS)$, by $(ii)$. 
\end{proof}

\bmhead{Acknowledgements}
The author wishes to thank Prof. Paul Deheuvels, who suggested this problem to him, and Prof. Zhan Shi for helpful discussions. The author also thanks Dr Brice Hannebicque for studying condition $\left( \Delta^{d_{\cH}} \right)$. 

\bmhead{Data availability}
No data is used.


\section*{Declarations}

\bmhead{Funding}
No funding.

\bmhead{Conflict of interest}
The author reports there is no conflict of interest to declare.

\bibliography{new}% common bib file

\end{document}